\documentclass[11pt,a4paper,leqno]{amsart}

\usepackage{fullpage,amsmath,amsthm,amsfonts,amscd,graphicx,
latexsym,amssymb,eucal,curves}
  \usepackage[all]{xy}
\newcommand{\ZZ}{\mathbb{Z}}

\newcommand{\CC}{\mathbb{C}}
\newcommand{\PP}{\mathbb{P}}
\newcommand{\NN}{\mathbb{N}}

\newcommand{\HH}{\mathbb{H}}

\newcommand{\LL}{\mathbb{L}}
\newcommand{\MM}{\mathbb{M}}

\newcommand{\WW}{\mathbb{W}}

\newcommand{\QQ}{\mathbb{Q}}
\newcommand{\RR}{\mathbb{R}}

\newcommand{\VV}{\mathbb{V}}
\newcommand{\cE}{\mathcal{E}}
\newcommand{\cF}{\mathcal{F}}
\newcommand{\cL}{\mathcal{L}}
\newcommand{\cO}{\mathcal{O}}

\newcommand{\Ker}{{\rm Ker}}

\newcommand{\JacS}{{\rm Jac}}

\newcommand{\Gal}{{\rm Gal}}

\newcommand{\Aff}{{\rm Aff}}

\newcommand{\Norm}{{\rm Norm}}
\newcommand{\tr}{{\rm tr}}

\newcommand{\MW}{{\rm MW}}

\newcommand{\SL}{{\rm SL}}
\newcommand{\GL}{{\rm GL}}
\newcommand{\diag}{{\rm diag}}

\newcommand{\Sing}{{\rm Sing}}

\newcommand{\prim}{{\rm prim}}
\newcommand{\sms}{\smallsetminus}

\newcommand{\ol}{\overline}

\newtheorem{Defi}{Definition}[section]

\newtheorem{Prop}[Defi]{Proposition}
\newtheorem{Lemma}[Defi]{Lemma}
\newtheorem{Cor}[Defi]{Corollary}
\newtheorem{Thm}[Defi]{Theorem}

\begin{document}
\title{Periodic points on Veech surfaces and 
the Mordell-Weil group over a Teichm\"uller curve}
\author{Martin M\"oller}

\begin{abstract}
Periodic points are points on Veech surfaces, whose
orbit under the group of affine diffeomorphisms is finite. 
We characterise those points as being torsion points
if the Veech surfaces is suitably mapped to its Jacobian
or an appropriate factor thereof. For a primitive
Veech surface in genus two we show that the only periodic
points are the Weierstra\ss\ points and the singularities.
\newline
Our main tool is the Hodge-theoretic characterisation of
Teichm\"uller curves. We deduce from it a finiteness result
for the Mordell-Weil group of the family of Jacobians
over a Teichm\"uller curve.
\end{abstract}
\date{\today}
\thanks{Supported by the 
DFG-Schwerpunkt ``Komplexe Mannigfaltigkeiten''}
\maketitle

\section*{Introduction}
Let $\Omega M_g$ denote the tautological bundle over the
moduli space of curves of genus $g$. Its points consist
of pairs $(X^0, \omega^0)$ of a Riemann surface plus a
holomorphic one-form on $X^0$. There is a natural $\GL^+_2(\RR)$-action
on $\Omega M_g$ (see \cite{Ve89} or \cite{McM03}). 
In the rare cases where the projection
of $\GL^+_2(\RR) \cdot (X^0, \omega^0)$ to $M_g$ is an
algebraic curve $C$, this curve is called a {\em Teichm\"uller curve}
and $(X^0, \omega^0)$ a {\em Veech surface}.
These surfaces are also characterised by the following
property:
\newline
Let $\Aff(X^0, \omega^0)$ denote the group of 
diffeomorphisms on $X^0$ that are affine outside the zeros of $\omega^0$
with respect to the charts determined
by integrating $\omega^0$. Then $(X^0, \omega^0)$ is a Veech surface
if and only if $\Aff(X^0, \omega^0)$ is 'as big as possible',
i.e.\ the matrix parts of these diffeomorphisms form
a lattice $\Gamma$ in $\SL_2(\RR)$.
\par
A {\em periodic point} on a Veech surface is a point whose
orbit under $\Aff(X^0, \omega^0)$ is finite. Examples of 
periodic points are the zeroes of $\omega^0$ and Weierstra\ss\ 
points if $g=2$ (see \cite{GuHuSc03}). 
\newline
The aim of this paper is to show:
\par
{\bf Theorem \ref{peristor}:} {\em The difference of
two periodic points is torsion on an $r$-dimensional
quotient of $\JacS(X^0)$, where $r= [\QQ(\tr(\Gamma)): \QQ]$. }
\par
We also show a converse
to this statement in terms of torsion sections of the family of 
Veech surfaces over an unramified cover of the Teichm\"uller 
curve. We say for short that periodic
points form a {\em torsion packet} (\cite{Co85}) on $X^0$.
\newline
As a consequence we also obtain a different proof of the
finiteness of the number of periodic points on a Teichm\"uller
curve that does not arise via a torus covering. This
was proved by Gutkin, Hubert and Schmidt (see \cite{GuHuSc03})
using the flat geometry induced by $\omega^0$.
\par
We apply the characterisation of periodic points via
torsion points to show:
\par
{\bf Theorem \ref{PerStratum2}/\ref{PerPtDec}:} 
{\em The only periodic points on a Veech surface
in genus $2$, which does not arise from a torus cover, are the
Weierstra\ss\ points and the zeroes of $\omega^0$. }
\par
Since 'large' torsion packets are 'rare' (see e.g.\ \cite{BoGr00}), 
we expect that for
large $g$ the strata of $\Omega M_g$ with 'many' (in particular 
$2g-2$ simple) zeroes should contain 'few' Teich\-m\"uller curves.
\newline
In fact for $g=2$ McMullen shows (\cite{McM04b} and \cite{McM04c}),
using the above characterisation of periodic points, that in $\Omega M_2(1,1)$
there is only one Teichm\"uller curve not coming from genus $1$. 
It is generated by the decagon. In contrast to that $\Omega M_2(2)$ 
contains infinitely many Teichm\"uller curves that are not obtained 
via torus coverings (see \cite{McM04a}).
\par
{\bf Contents.} \newline
In Section \ref{dict} we start by recalling the language
of translation surfaces in which periodic points were
first studied in \cite{GuHuSc03}. We translate this into the
following setting:
\newline
To a Veech surface we may associate (via its $\GL^+_2(\RR)$-orbit) 
a Teichm\"uller curve $C \to M_g$. 
A Teichm\"uller curve defines (by pulling back the universal 
family over a finite cover of $M_g$) a family of curves 
(or a {\em fibred surface})
$f: X_1 \to C_1$, after passing to some finite 
unramified covering $C_1$ of $C$. Finally each fibre
of $f$ is a Veech surface $X^0$ and the $1$-form
$\omega^0$ can be recovered from the variation of
Hodge structures of $f$.
\newline
We show that periodic points correspond via this
dictionary to sections of $f$.
\par
In Section \ref{PrimIrr} we recall some material from 
\cite{Mo04} in order to define the $r$-dimensional factor $A_1/C_1$ of
$\JacS(X_1/C_1)$ we referred to above. 
\newline
We will call a Veech surface {\em geometrically primitive} if it 
does not arise via coverings from smaller genus. We will call
a Veech surface {\em algebraically primitive} if the relative Jacobian
of an associated fibred surface is irreducible. The indefinite
article here and whenever we talk of fibred surfaces refers
only to the possibility of passing to an unramified covering.
\newline
We give an example that (unlike in genus two) these notions
do not coincide in higher genus. In Theorem \ref{uniqueroot}
we show that nevertheless a Veech surface stems from
a unique geometrically primitive Veech surface.
\par 
The {\em Mordell-Weil group} $\MW(A/C)$ of a family of abelian varieties
$g: A \to C$ is the group of rational sections of $g$ or
equivalently the group of $\CC(C)$-valued points of $A$. 
In Section \ref{MWgroup} we
study the Mordell-Weil group for the factor $A_1/C_1$ of $\JacS(X_1/C_1)$
for a Teichm\"uller curve $C \to M_g$:
\newline
{\bf Theorem \ref{MordellWeil}:} {\em For each Teichm\"uller curve and any
given unramified covering $C_1 \to C$ the group
$\MW(A_1/C_1)$ is finite. }
\par
We apply this to the characterisation
of periodic points via torsion sections (Thm.\ \ref{peristor}
and the converse Prop.\ \ref{converse}).
\par
The last two sections contain a degeneration
argument similar to the one in \cite{McM04c} for 
the explicit analysis of periodic points in genus two. 
While in loc.\ cit.\ ratios of sines appear, we need for our
purposes that only finitely many ratios 
$\tan(a)/\tan(b)$
for $a,b \in \pi\QQ$ lie in a quadratic number field. 
We give an explicit list of them, that is used in the last section.
\par
{\bf Some notation.} 
\newline
Riemann surfaces are usually
denoted by $X^0$, $Y^0$, etc., while $X$, $X_1$, $Y$, $Y_1$, etc.\ will
be used for families of Riemann surfaces over base curves $C$, $C_1$, etc.,
whose fibres over $0 \in C(\CC)$ are $X^0$, $Y^0$, etc.
\newline
Overlines denote completions of curves or fibred surfaces.
\par
{\bf Acknowledgements.}
\newline
The author thanks P.~Hubert, whose questions motivated this paper, and
C.~McMullen for useful discussion.
\par

\section{A dictionary} \label{dict}

{\bf Translation surfaces. } \newline
A {\em translation surface} is a compact Riemann surface $X^0$ with
an atlas $\{U_i, i\in I\}$, covering $X^0$ except for finitely many 
points (called {\em singularities} $\Sing(X^0)$) such that 
the transition functions between the charts are translations. 
A holomorphic one-form $\omega^0$ on a Riemann surface 
$X_0$ induces the structure of a translation
surface where $\Sing(X^0)$ can be any finite set containing
the zeroes of $\omega^0$.
A {\em translation covering} between two translation surfaces
$(X^0,\Sing(X^0))$ and $(Y^0,\Sing(Y^0))$ is a covering
$\varphi:X^0 \to Y^0$ such that the charts of $X^0$ are pulled back
from charts of $Y^0$ and such that $\varphi^{-1}(\Sing(Y^0)) =
\Sing(X^0)$. These coverings are sometimes called {\em balanced}
translation coverings. We deal exclusively with them.
\newline
If the translation structure on $Y^0$ is induced from $\omega_Y^0$ 
and $\varphi$ is a translation covering then the
translation structure on $X^0$ is induced from $\omega^0 = 
\varphi^* \omega_Y^0$. 
\newline
Finally note that a (possibly ramified) covering $\varphi: X^0 \to
Y^0$ plus differentials on $X^0$ and $Y^0$ satisfying
$\omega^0 = \varphi^* \omega^0_Y$ induces a translation covering with
$\Sing(Y^0)$ the zeros of $\omega_Y^0$ and $\Sing(X^0)$ their
preimages. This set will in general properly contain the zeroes of $\omega^0$.
We will abbreviate a translation covering by
$\pi: (X^0, \omega_X^0) \to (Y^0,  \omega_Y^0)$ or
$\pi: (X^0, \Sing(X^0)) \to (Y^0, \Sing(Y^0)$. 
\par
{\bf Affine diffeomorphisms, affine group.}
\newline
On a translation surface we may consider diffeomorphisms
that are orientation-preserving, 
affine with respect to the coordinate charts and that
permute the singularities.
They form the {\em group of affine diffeomorphisms} denoted by
$\Aff(X^0, \Sing(X^0))$. The matrix part 
of such an affine diffeomorphism is well-defined and
this yields a map $D$ to 
a subgroup  $\Gamma$ in $\SL_2(\RR)$. We call $\Gamma$ 
the {\em affine group} of $(X^0, \Sing(X^0))$.
If $\Sing(X^0)$ consists of the zeroes of $\omega^0$ we sometimes
write $\SL(X^0,\omega^0)$ for $\Gamma$.  
The kernel of $D$ consists of conformal automorphisms of $X^0$
preserving $\omega^0$ and is hence finite.
\par
\begin{Defi}
A point $P$ on a translation surface is called {\em periodic} 
if its orbit under the group $\Aff(X^0,\Sing(X^0))$ is finite.
\end{Defi}
\par
{\bf Veech surfaces, Teichm\"uller curves.}
\newline
If $\Gamma \subset \SL_2(\RR)$ is a lattice, the
translation surface $(X^0,\Sing(X^0))$ is called
a {\em Veech surface.} If $\Sing(X^0)$ consists of the zeroes of 
$\omega^0$, this is the case if and only if the $\GL^+_2(\RR)$-orbit
of $(X^0,\omega^0)$ in $\Omega M_g$ projects to 
an algebraic curve $\HH/\Gamma = C \to M_g$, which is called a 
{\em Teichm\"uller curve.} The map $C \to M_g$ is injective
up to finitely many normal crossings. By abuse of notation
we will also call Teichm\"uller curve a map $C_1 \to M_g$, 
which is the composition of an unramified cover $C_1 \to C$
composed by a Teichm\"uller curve $C \to M_g$ in the above
sense.
\par
If we restrict to the quotient $C_1 = \HH/\Gamma_1$ for
a sufficiently small subgroup $\Gamma_1 \subset \Gamma$
of finite index, $C_1$ will map to a finite cover of 
$M_g$ over which the universal family exists.
If we pull back this universal family to $C_1$ we
obtain a {\em fibred surface $f:X_1 \to C_1$} associated
with the Teichm\"uller curve. We will also need 
a smooth semistable model $f:\ol{X_1} \to \ol{C_1}$ 
over the completion of $C_1$. 
\newline
A Veech surface $(X^0, \Sing(X^0))$ 
is called {\em square-tiled} (or {\em arises as torus cover}, 
or {\em origami}) 
if it admits a translation covering to a torus with one singular
point.
\par
{\em From now on we will exclusively deal with Teichm\"uller
curves and translation surfaces that are Veech surfaces.}
\par
Note that 'fibred surface' refers to an object of
complex dimension two. It contains as one of its
fibres the translation surface $X^0$, an object
of real dimension two.
\par
\begin{Lemma} \label{perfib}
A point $P$ on a Veech surface $(X^0, \omega^0)$ is periodic if and only if
there is an (algebraic) section of some
fibred surface $f:X_1 \to C_1$ associated
to the Teichm\"uller curve, which passes through $P$ on the
fibre $X^0$ of $f$. 
\end{Lemma}
\par
{\bf Proof:}
A section $s$ of $f$ over $C_1 = \HH/\Gamma_1$ hits $X^0$
in one point $P$. The $\Gamma$-orbit of $P$ consists
of at most $[\Gamma:\Gamma_1]$ points, hence is finite
by the choice of $\Gamma_1$.
\newline
Conversely given a periodic point $P$ on $X^0$ we may take a double
cover $\pi: Y^0 \to X^0$ branched at $P$ and some zero of $\omega^0$.
The translation surface $(Y^0, \omega_Y^0 := \pi^* \omega^0)$ 
is still a Veech surface: 
Indeed let $\Sing(X^0) = Z(\omega^0) \cup \{P\}$. 
The affine group of $(X^0, \Sing(X^0))$, 
is of finite index in $\Gamma$ by the definition of a periodic point.
Now $\pi$ defines a translation covering $(Y^0, \pi^{-1}(\Sing(X^0)))
\to (X^0, \Sing(X^0))$ and one can apply \cite{GuJu00} Thm.\ 4.9
to show the Veech property of $(Y^0, \omega_Y^0)$.
\newline  
Over some $C_1 = \HH/\Gamma_1$ for a subgroup $\Gamma_1$ of 
finite index in $\Gamma$ we have a covering $\pi: Y_1 \to X_1$
of fibred surfaces over $C_1$, such that the original $\pi$ is
the fibre over some point $0 \in C_1(\CC)$. By construction
of $Y_1$ as $\SL_2(\RR)$-orbit of $(Y^0, \omega_Y^0)$ the differential
$\omega_Y^0$ extends to a section $\omega_{Y_1}$ of the relative canonical
sheaf $\omega_{\ol{Y_1}/\ol{C_1}}$. Again by definition of
the $\SL_2(\RR)$-action the multiplicities of the zeros of $\omega_{Y_1}$
remain constant over $C_1$. Hence passing to a subgroup of
finite index in $\Gamma_1$ (we nevertheless keep the notation)
we may assume that the zeros of $\omega_{Y_1}$ define sections $s_i$
of $f_Y: Y_1 \to C_1$. 
\newline
The images of $s_i$ under $\pi$ are sections of $f$. One of them passes
through $P$, as $\pi$ is ramified over $P$.
\hfill $\Box$
\par

\section{Algebraic and geometric primitivity} \label{PrimIrr}

Let $K = \QQ(\tr(\Gamma))$ denote the trace field of the
affine group of a Veech surface $(X^0, \omega^0)$. It remains unchanged
if we replace $\Gamma$ by a subgroup of finite
index. Let $r := [K:\QQ]$. We recall from \cite{Mo04} 
the decomposition of the variation of Hodge structures (VHS)
over a Teichm\"uller curve generated by $(X^0, \omega^0)$:
\newline
Let $\VV = R^1 f_* \ZZ$. In \cite{Mo04} Prop.\ 2.3 we
have shown that there is a decomposition as polarized VHS 
$$ \VV_\QQ = \WW \oplus \MM, \quad \WW_K = \bigoplus_{\sigma \in \Gal(K/\QQ)}
\LL^\sigma
\eqno(1)$$
Here $\LL^\sigma$ are rank two local systems over $K$ and
none of the irreducible factors of $\MM_\CC$ is isomorphic 
to any of the $\LL_\CC^\sigma$.
This yields (see loc.\ cit.\ Thm.\ 2.5) a decomposition of the Jacobian 
$$\JacS(X_1/C_1) \sim 
A_1 \times B_1 \eqno(2)$$ 
up to isogeny, where $A_1$ has dimension $r$
and real multiplication by $K$. 
\newline
Recall furthermore that the graded quotients of a VHS together
with the Gauss-Manin connection form a Higgs bundle
$({\cE}, \Theta)$. The summands $\LL^\sigma$ of $\VV_K$
give rank-two sub-Higgs bundles $({\cL}^\sigma \oplus
({\cL}^\sigma)^{-1}), \tau^\sigma)$, where $S = \ol{C} \sms C$ and
$$\tau^\sigma: {\cL}^\sigma \to ({\cL}^\sigma )^{-1} \otimes 
\Omega^1_{\ol{C}}(\log S).$$
The subbundle ${\cL}^{\rm id}$ of 
$f_* \omega_{\ol{X_1}/\ol{C_1}}$ is distinguished by the property that
its restriction to the fibre $X^0$ gives $\CC \cdot \omega^0$.
\newline
Teichm\"uller curves are characterised (see \cite{Mo04} Thm.\ 5.3)
by a decomposition of the VHS as above plus the property that
$\LL^{\rm id}$ is {\em maximal Higgs}, i.e.\ that $\tau^{\rm id}$
is an isomorphism. We need only two properties
of this notion here: It is stable under replacing $C_1$ by a finite
unramified cover and the VHS over a Teichm\"uller curve 
has precisely one rank-two subbundle that is maximal Higgs.
\newline
The last property is stated explicitely in \cite{Mo04} Lemma 3.1 
for the $\LL^{\sigma}$. But any two maximal Higgs subbundles become
isomorphic after replacing $C_1$ by a finite \'etale cover
and none of the irreducible summands of $\MM_\CC$ is
isomorphic to any of the $\LL^\sigma_\CC$.
\par
\begin{Defi}
A Teichm\"uller curve is called {\em algebraically primitive} if its 
relative Jacobian $\JacS(X_1/C_1)$ is irreducible as abelian
variety.
\end{Defi}
\par
We will also say that a Veech surface is algebraically primitive,
if the corresponding Teich\-m\"uller curve is algebraically primitive.
\begin{Lemma}
A Teichm\"uller curve is algebraically primitive if and only if
$r=g$. 
\end{Lemma}
\par
{\bf Proof:} 
We only have to show that $\WW$ is irreducible over $\QQ$.
This follows immedately from the irreducibility of the $\LL^\sigma$
and the fact that for $\sigma \neq \tau$ the local systems
$\LL^\sigma$ and $\LL^\tau$ are not isomorphic (see 
\cite{Mo04} Lemma 2.2). \hfill $\Box$
\par
Note that irreducibility of $\JacS(X_1/C_1)$
does not depend on replacing
$C_1$ by unramified covers. 
\newline
Further note that irreducibility of $\JacS(X_1/C_1)$ does 
not exclude that special fibres of $f: X \to C$ may have reducible
Jacobians. 
\par
There is also a natural geometric notion of primitivity for
translation surfaces and one for Teichm\"uller curves without explicitly
referring to any differental. We show that
these two geometric definitions coincide.
\par
\begin{Defi} A translation surface $(X^0,\omega_X^0)$
is {\em geometrically imprimitive} if there exists a
translation surface $(Y^0, \omega_Y^0)$ of smaller genus
and a covering $\pi :X^0 \to Y^0$ such that
$\pi^* \omega_Y^0 = \omega_X^0$.
\end{Defi}
\par
\begin{Defi}
A Teichm\"uller curve $C \to M_g$ is {\em geometrically imprimitive} if
there is an unramified cover $C_1 \to C$, a
fibred surface $f_Y: Y_1 \to C_1$ coming from
a Teichm\"uller curve $C_1 \to M_h$ with $h < g$ and
a (possibly ramified) covering $\pi: X_1 \to Y_1$ over $C_1$.
It is called {\em geometrically primitive} otherwise.
\end{Defi}
\par
\begin{Prop} \label{PrimEquiv}
A Teichm\"uller curve is geometrically primitive if and only if a
corresponding Veech surface is geometrically primitive.
\end{Prop}
\par
{\bf Proof:} 
Suppose that $C \to M_g$ is geometrically imprimitive and 
let $\pi: X_1 \to Y_1$
be the covering of fibred surfaces associated with the
Teichm\"uller curves.  We restrict 
the covering $\pi$ to some fibre $X_1^0 \to Y_1^0$
over $0 \in C_1(\CC)$. 
We have to show that the differentials $\omega_{X_1}^0$ on
$X_1^0$ and $\omega_{Y_1}^0$ on $Y_1^0$ that
generate the Teichm\"uller curves define 
a translation covering. 
We let $\Sing(Y_1^0)$ 
be the zeroes of $\omega_{Y_1}^0$ and $\Sing(X_1^0)$
be their preimage via $\pi$. It suffices
to see that (up to a multiplicative constant)
we have $\pi^* \omega_{Y_1}^0 = \omega_{X_1}^0$.
By definition $\omega_{X_1}^0$ and $\omega_{Y_1}^0$ 
are obtained from ${\cL}^{\rm id}_{X_1}$ resp.\
 ${\cL}^{\rm id}_{Y_1}$ by restriction to $0$. 
By the properties of maximal Higgs subbundles listed
above, we have ${\cL}^{\rm id}_{X_1} = \pi^* 
{\cL}^{\rm id}_{Y_1}$. This  completes the 'if'-part.
\par
Conversely a translation covering defines two
Teichm\"uller curves with commensurable affine groups
(\cite{Vo96} Thm.\ 5.4 or \cite{GuJu00} Thm.\ 4.9).
We pass to a subgroup $\Gamma_1$ contained in both with
finite index and small enough to have universal
families. Then the $\SL_2(\RR)$-images of the
translation covering patch together to a
covering of fibred surfaces over $\HH/\Gamma_1$.
\hfill $\Box$
\par
{\bf Examples:} i) The two notions of primitivity coincide
for genus $2$ and will therefore be abbreviated just by
'primitive': If $r=1$ the factor $A$ of the
Jacobian is one-dimensional, hence a family of
elliptic curves and the projection onto
this factor implies geometric imprimitivity (\cite{GuJu00} Thm.\ 5.9). 
\newline
ii) The Riemann surface $y^2 = x^7-1$ with $\omega_0 =dx/y$
is a Veech surface (see \cite{Ve89}) and $\Gamma = \Delta(2,7,\infty)$
has a trace field of degree $3$ over $\QQ$. The corresponding
Teichm\"uller curve is algebraically primitive, hence 
geometrically primitive. 
\newline
iii) For genus $g=3$ the notions of geometric and
algebraic primitivity no longer coincide: The Riemann surface 
$y^{12} = x^3(x-1)^4$ with 
$\omega^0 = ydx/[x(x-1)]$ is studied in \cite{HuSc01} Thm.~5.
It is obtained from unfolding the billiard in the triangle
with angles $3\pi/12$, $4\pi/12$ and $5\pi/12$. It is
a Veech surface with $\Gamma=\Delta(6,\infty,\infty)$.
The trace field is of degree only $2$ over $\QQ$. Hence the 
Teichm\"uller curve is not algebraically primitive. But it is 
geometrically primitive, 
as remarked in loc.\ cit.: Since $r=2$ it cannot arise as a torus cover.
If it arose as a genus two cover, this cover would have to be
unramified by Riemann-Hurwitz. Hence $\omega^0$ would
have zeros of order at most two. But $\omega^0$ has indeed
a zero of order four.
\par
In \cite{HuSc01} the authors analyse the properties
that are preserved if one goes up and down a tree
of translation coverings. We show that the situation
is simple, if the singularities of the translation
surfaces are chosen suitably, i.e.\ each tree has
a root. 
\par
\begin{Thm} \label{uniqueroot}
A translation surface $(X^0,\omega^0)$ is obtained as a translation
covering from a geometrically primitive translation surface
$(X^0_\prim, \omega^0_\prim)$. If the genus $X^0_\prim$ is greater than one,
this primitive surface is unique.
\newline
Moreover, the construction of the primitve surface is equivariant with
respect to the $\SL_2(\RR)$-action, if the genus of $X^0_\prim$ is greater
than one. 
\newline
If $(X^0,\omega^0)$ generates a Teichm\"uller curve,
then so does $(X^0_\prim, \omega^0_\prim)$. 
In this case $\pi: X^0 \to X^0_\prim$ is branched over periodic points.
\end{Thm}
\par
{\bf Proof:} We first prove existence and uniquenes of the 
primitive surface
and drop the superscripts that indicate a special fibre.
\newline
Let $\JacS(X) \to A^{\rm max}_X$ be the maximal abelian quotient such that
$\omega_X$ pulls back from $A^{\rm max}_X$. 
Equivalently $A^{\rm max}_X$ is the quotient of $\JacS(X)$
by all connected abelian subvarieties $B'$ such that
the pullback of $\omega_X$ to $B'$ vanishes.  
If $\dim A^{\rm max}_X=1$ this is the primitive surface, we are done.
\newline
Else embed $X$ into its Jacobian and consider the normalization $X_A$
of the images of $X \to \JacS(X) \to A_X^{\rm max} \to A$ for
all isogenies $A_X^{\rm max} \to A$. Since $X_A$ generates $A$
as a group, we have $g(X_A) \geq \dim A \geq 2$. The curves $X_A$ form 
an inductive system which eventually stabilizes since the genera
are bounded below. We let $X_\prim$ be the limit of the $X_A$
and claim that it has the required properties.
\newline
First, by construction of $A$ and since $X_\prim$ generates
$A$ as a group, there is a differential $\omega_\prim$
on $X_\prim$ such that $\varphi_X^* \omega_\prim = \omega_X$.
Hence $\varphi_{X}: X \to X_\prim$ defines a translation cover
once the singularities are suitably chosen.
\newline
Second, $(X_\prim,\omega_\prim)$ is indeed primitve:
Suppose there is a covering $X_\prim \to Y$ with 
$\pi^* \omega_Y = \omega_\prim$ for some differential on $Y$. Let $\JacS(\pi)$ be
the induced morphism on the Jacobians, commuting with $\pi$
and suitable embeddings of the curves into their Jacobians.
Consinder the following diagram:
$$ 
\xymatrix{ X_\prim \ar[d]^{\pi} \ar@{-->}[drrr]^{\cong} 
\ar[r] & \JacS(X) \ar[rrr] \ar[dd]^>>>>>>{\JacS(\pi)}  
&&& A_X^{\rm max} \ar@{-->}[dd]^>>>>>>{\ol{\JacS(\pi)}}_>>>>>>{\exists} \\
Y  \ar[dr] \ar[rrr] &&& Z \ar[dr] \\
& \JacS(Y) \ar[rrr] &&& A_Y^{\rm max} \\
}$$
By definition $\omega_\prim$ vanishes on $K :=
\Ker(\JacS(X) \to A^{\rm max}_X)$. Since $\pi^* \omega_{Y} = \omega_\prim$
the differential $\omega_{Y}$ vanishes on $\JacS(\pi)(K)$.
This means that $\JacS(\pi)(K)$ is in the kernel of
$\JacS(Y)\to A(Y)$ and hence $\JacS(\pi)$ descends
to a homomorphism on the quotients $\ol{\JacS(\pi)}: A^{\rm max}_X \to A^{\rm
max}_Y$. This is an isogeny by construction of $A^{\rm max}_X$.
The map $X_\prim \to Z$, which is defined as the normalization of the
image in $A^{\rm max}_Y$, is an isomorphism by construction and so is $\pi$.
\newline
Third, for the uniqueness we have to show that
for a translation covering $\pi:(X,\omega_X) \to (Y,\omega_Y)$ 
there is a morphism $\pi_\prim: X_\prim \to Y_\prim$.
As in the second step we have an induced map 
$\ol{\JacS(\pi)}: A_X^{\rm max} \to A_Y^{\rm max}$.
The curve $X_\prim$ was obtained as the normalization of
the image of $X$ is some quotient $q_X: A_X^{\rm max} \to A_X$.
It hence maps to the normalization $Z$ of the image of $Y$ in
$A_X^{\rm max} /\langle\Ker(q_X),\Ker(\ol{\JacS(\pi))}\rangle$.
Since $Z$ maps to $Y_\prim$ by construction we are done.
\par
For the $\SL_2(\RR)$-equivariance let $(X^1,\omega^1) = A\cdot(X^0,\omega^0)$
for some $A \in \SL_2(\RR)$. Primitivity implies the existence of a
translation coverings $\pi: A\cdot X^0_\prim \to X^1_\prim$ and
$\pi': A^{-1}\cdot X^1_\prim \to X^0_\prim$. Hence either both primitive
curves have genus $1$ or both have bigger genus greater than one
and we are done by uniqueness.
\par
In case of a torus cover the statement that
$(X^0_\prim,\omega^0_\prim)$ generates a Teichm\"uller curve
is trivial. In the other cases the previous argument implies
that the affine group of $(X^0_\prim, \omega^0_\prim)$ contains
the one of $(X^0,\omega^0)$.
\par
If both translation surfaces generate Teichm\"uller curves
there is a subgroup of finite index of $\Aff(X^0, \omega^0)$
that descends to $X^0_\prim$. This group has to fix branch points
and hence the whole group can generate only finite orbits of
branch points.
\hfill $\Box$
\par

\section{The Mordell-Weil group} \label{MWgroup}

\begin{Thm} \label{MordellWeil}
Let $f: X_1 \to C_1$ be 
a fibred surface associated with a Teichm\"uller curve. 
Then the Mordell-Weil group of $A_1/C_1$ is finite. 
Here $A_1$ is the factor of $\JacS(X_1/C_1)$
with real multiplication by $K$.
\newline
In particular if the Teichm\"uller curve
is algebraically primitive then  $\MW(\JacS(X_1/C_1))$ is finite.
\end{Thm}
\par
The finiteness of the Mordell-Weil group is
invariant under isogenies. Thus there is no need to specify
$A_1$ in its isogeny class. In particular we may suppose 
that the $\QQ$-local system $\WW$ comes from a $\ZZ$-local system
$\WW_\ZZ$.
\newline
Furthermore the statement of the theorem becomes stronger the
smaller the subgroup $\Gamma_1$ with $C_1 = \HH/\Gamma_1$ is.
Thus we may replace $C_1$ by an unramified cover and suppose
that we have unipotent monodromies. To simplify
notation we will call this cover $C$, which should
not be confused with the notation for the original
Teichm\"uller curve. 
\newline
Let $\ol{g}:\ol{A} \to \ol{C}$ be an extension of $g: A \to C$ to
a semiabelian scheme. A unique such extension exists
due to the unipotent monodromies. We denote by $H^0(C, {\cO}_{C}(A/C))$
(resp.\ by $H^0(\ol{C},{\cO}_{\ol{C}}^{\rm an}(\ol{A}))$)
the group of algebraic sections of $A/C$ (resp.\ the
group of analytic sections of $\ol{A}/\ol{C}$).
\newline
Two remarks: First, the analytic sections of  $\ol{A}/\ol{C}$ are 
necessarily algebraic. Nevertheless we have to use the analytic
category, because we want to use uniformization in the sequel.
Second, by properness of $g$ any rational section of $g$ extends
to the whole curve $C$, hence $\MW(A/C) = H^0(C, {\cO}_{C}(A/C))$. 
\par
\begin{Lemma} 
The restriction map
$$r:H^0(\ol{C},{\cO}_{\ol{C}}^{\rm an}(\ol{A})) \to
H^0(C, {\cO}_{C}(A/C)) $$
is injective with finite cokernel.
\end{Lemma}
\par
{\bf Proof:} The proof is from Prop.\ 6.8 in \cite{Sa93}. 
We reproduce a sketch for convenience of the reader:
\newline
There exists a group scheme $N \to \ol{C}$ (called the N\'eron model 
of $\ol{A}/\ol{C}$) containing $\ol{A}/\ol{C}$
as a connected component, with the following property: For smooth 
$Y \to \ol{C}$ a rational map $Y --> N$ over $\ol{C}$ 
extends to a morphism $Y \to N$.
\newline
We apply this property to sections $s: C \to A$. They extend
to rational maps from $\ol{C}$ a priori to a projective completion 
of $N$, but in fact to a morphism $\ol{C} \to N$. Hence
$$ H^0(C, {\cO}_{C}(A/C)) \to H^0(\ol{C},{\cO}_{\ol{C}}(N/\ol{C}))$$
is an isomorphism. This shows that the cokernel of $r$ consists of 
sections of the finite group scheme $N/\ol{A}$ and is hence finite. 
\hfill $\Box$
\par
If the local system $\WW$ on $C$ carries a polarized VHS of weight $m$ 
(in our case $m=1$)
then $H^i(\ol{C},j_* \WW_\CC)$ is known to carry a Hodge structure of
weight $m+i$ (see \cite{Zu79}). Here $j: C \to \ol{C}$ is
the inclusion. Indeed let $\Omega^\bullet(\WW)_{(2)}$ 
denote the deRham complex with $L_2$-growth
conditions at the punctures. Zucker (see \cite{Zu79}) shows 
(extending Deligne's
results to the non-compact case) that we may identify 
$H^i(\ol{C},j_* \WW_\CC)$ with the hypercohomology groups
${\bf H}^i(\ol{C}, \Omega^\bullet(\WW)_{(2)})$. Then the Hodge structure comes 
from the filtration
on $\Omega^\bullet(\WW)_{(2)}$ induced by the Hodge filtration on $\WW$.
\newline
But in fact there is another complex more easily accessible
and quasi-isomorphic to $\Omega^\bullet(\WW)_{(2)}$ (\cite{Zu79} Prop.\ 9.1).
\newline
We describe how to calculate the Hodge structure on $H^i(\ol{C},j_* \WW_\CC)$
in our situation: The Hodge filtration on $\WW$ is
$$ 0={\cF}^2 \subset {\cF}^1 = g_* \omega_{\ol{A}/\ol{C}}
\subset {\cF}^0 = (\WW \otimes {\cO}_C)_{\rm ext}, $$
where the subscript denotes the Deligne extension to $\ol{C}$. 
As graded pieces we have 
$$ {\cE}^{1,0} = {\cF}^1/{\cF}^2 =  g_* \omega_{\ol{A}/\ol{C}},
\quad {\cE}^{0,1} = {\cF}^0/{\cF}^1 =  R^1 g_* {\cO}_{\ol{A}}. $$
\par
Combining \cite{Zu79} Thm.\ 7.13 and Prop.\ 9.1 (see also the
restatement after Lemma 12.14 in loc.\ cit.) we conclude for $p \in \{0,1,2\}$ 
$$ H^1(\ol{C},j_* \WW_\CC)^{(p,2-p)} = {\bf H}^1(\ol{C}, {\cE}^{p,1-p} 
\longrightarrow 
{\cE}^{p-1,2-p} \otimes \Omega^1_{\ol{C}}(\log S)) \eqno(3)$$
where the mapping in the complex on the right is the graded quotient
of the Gauss-Manin connection (equivalently: the 
Kodaira-Spencer mapping) and $S = \ol{C} \sms C$. 
\par
\par
{\bf Proof of the theorem:} 
The uniformization of $\ol{A}/\ol{C}$ yields a short exact sequence
$$ 0 \to j_* \WW_\ZZ \to {\cE}^{0,1} \to {\cO}_{\ol{C}}^{\rm an}(\ol{A}) 
\to 0.$$
We take cohomology and note that $H^0(\ol{C}, {\cE}^{0,1})$ vanishes as
$\ol{A}/\ol{C}$ has no fixed part. Hence
$$
\begin{array}{lcl}
H^0(\ol{C}, {\cO}_{\ol{C}}^{\rm an}(\ol{A})) & = & \Ker 
(H^1(\ol{C},j_* \WW_\ZZ) \to H^1(\ol{C}, {\cE}^{0,1})) \\
& = & H^1(\ol{C},
j_* \WW_\ZZ) \cap (\Ker 
(H^1(\ol{C}, j_* \WW_\CC) \to H^1(\ol{C}, {\cE}^{0,1}))) \\
& = & H^1(\ol{C},
j_* \WW_\ZZ) \cap H^1(\ol{C}, j_* \WW_\CC)^{1,1}. 
\end{array}
$$
By $(1)$ and $(3)$ we deduce
$$ H^1(\ol{C}, j_* \WW_\CC)^{1,1} = \oplus_{\sigma_i \in \Gal(K/\QQ)}
{\bf H}^1(\ol{C}, {\cL}_i \to ({\cL}_i)^{-1}  
\otimes \Omega^1_{\ol{C}}(\log S)),$$
where ${\cL}_i$ is the $(1,0)$-part of
$(\LL^{\sigma_i} \otimes {\cO}_C)_{{\rm ext}}$. 
As the Kodaira-Spencer map for $\sigma_1 := {\rm id}$
is an isomorphism (this is the definition of 'maximal Higgs'), 
the first summand vanishes. But the action of $K$
permutes the summands transitively and hence 
$H^1(\ol{C}, j_* \WW_\ZZ) \cap H^1(\ol{C}, j_* \WW_\CC)^{1,1} = 0$. 
\hfill $\Box$ 
\par
\begin{Thm} \label{peristor}
Let $\varphi: X^0 \to \JacS (X^0) \to A^0$ be the embedding
of a Veech surface into its Jacobian (normalized such 
that one of the zeros of $\omega^0$ maps to $0$)
composed by the projection to the factor $A^0$.
\newline
The periodic points on a Veech surface map via $\varphi$ to 
torsion points on $A^0$.
In particular there is only a finite number of periodic
points on a Veech surface if $r>1$, i.e.\
if the surface is not square-tiled.
\end{Thm}
\par
The finiteness result was obtained by Gutkin, Hubert and Schmidt by
entirely different methods in \cite{GuHuSc03}.
\par
{\bf Proof:}
A periodic point of $(X^0, \omega^0)$ gives a section of 
some fibred surface $X_1 \to C_1$ by Lemma \ref{perfib} and via $\varphi$
a section of $A_{C_1} \to C_1$. This section has finite order
by Thm.\ \ref{MordellWeil} thus proving the first statement. 
\newline
By Thm.\ 5.1 in \cite{Mo04} the family of abelian varieties $A/C$
and also the section are defined over some number field.
We fix some fibre of $f:X \to C$ defined over some number field,
say our original $X^0$. If $r>1$ the image $\varphi(X^0)$
in $A^0$ is a curve, which generates an abelian variety of 
dimension $r$. Hence
it cannot be (a translate of) an abelian subvariety. 
\newline
In this situation the Manin-Mumford conjecture says that
the $\ol{\QQ}$-rational points of $\varphi(X^0)$ have finite intersection
with $A^0_{\rm tors}$. We have seen that all periodic points are contained
in this intersection.
\newline
Proofs of the Manin-Mumford conjecture were obtained in different
generality by e.g.\ Raynaud, Hindry, Vojta, Buium, 
Hrushovski, McQuillan. For what we need here the proof
in \cite{PiRo02} is sufficient and maybe the most easily
accessible.
\hfill $\Box$
\par
\begin{Cor}
If $(X^0,\omega^0)$ is a Veech surface that generates an
algebraically primitive Teichm\"uller curve, then all periodic
points form a {\em torsion packet}, i.e.\ for two periodic 
points $P$, $Q$ the difference $P - Q$ is torsion
(as a divisor class).
\end{Cor}
\par
There is a converse to the above theorem, if we
look for torsion sections instead of looking fibrewise:
\par
\begin{Prop} \label{converse}
Let $\varphi_{C}: X \to \JacS(X/C) \to A$ be the
family of maps considered in the previous theorem
and let $\varphi_{C_1}: X_1 \to A_1$ be the map obtained by
an unramified base change $C_1 \to C$.
\newline
Periodic points on $X$ are precisely the preimages 
via $\varphi_{C_1}$ of sections of
$A_1 / C_1$ for all unramified coverings $C_1 \to C$.
\end{Prop}
\par
{\bf Proof:} We may choose $A$ as in the
proof of Thm.\ \ref{uniqueroot}, as the statement of the
proposition is invariant under isogenies.
Sections of $\varphi_{C_1}(X)$ extend to
sections of its normalisation, which was called
$f_\prim: X_\prim \to C_1$. Hence they give periodic
points on each fibre of $X_\prim$
by the criterion for periodic points given in Lemma \ref{perfib}.
Finally periodic points remain periodic under passing
to a translation cover of Veech surfaces. 
\hfill $\Box$
\par
\begin{Cor}
If $r>1$ there is a universal bound depending only on $g$ 
for the number of periodic points on a Veech surface
of genus $g$.
\end{Cor}
\par
In fact Buium (\cite{Bu94}) gives a bound on the
number of torsion sections of a family of curves
in a family of abelian varieties. This bound only
depends on the genus of the curve and the dimension
of the abelian variety, but it grows very fast with $g$.  
\par

\section{Ratios of tangents} \label{tans}

For the next section we will need:
\par
\begin{Thm}
For each $d>0$, there are only finitely many pairs of
rational numbers $0 < \alpha < \beta < 1/2$ such that
$$\mu = \tan(\pi\beta)/\tan(\pi\alpha)$$
is an algebraic number of degree $d$ over $\QQ$.
\end{Thm}
\par
{\bf Proof:} This follows from \cite{McM04c} Thm.\ 2.1 and
the addition formula
$$ \frac{\tan(\frac{x+y}{2})}{\tan(\frac{x-y}{2})} =
\frac{\frac{\sin(x)}{\sin(y)} + 1}{\frac{\sin(x)}{\sin(y)} - 1}.$$
\hfill $\Box$
\par
For $d=2$ a list of these quotients can be obtained easily
from table $3$ in \cite{McM04c}. For later use we list those
quotients, which are non-units. By Galois conjugation in the
cyclotomic field containing $\mu$ we may furthermore suppose $\alpha = 1/s$
for some $s \in \NN$.
\par
$$
\begin{array}{ccccc}
\alpha & \beta  & \mu & \text{Trace} 
& \text{Norm} \\ 
1/10 & 1/5 & \sqrt{5} & 0 & -5 \\
1/10 & 2/5 &  5 + 2 \sqrt{5} & 10 & 5 \\
1/5  & 3/10 & 1 + 2\sqrt{1/5}& 2 & 1/5  \\ 
1/12 & 1/6 &  (3 + 2\sqrt{3})/3 & 2 & -1/3 \\
1/12 & 1/3 &  3 + 2\sqrt{3} & 6 &  - 3 \\
1/6 &  1/4 & \sqrt{3} & 0 & -3  \\
1/6 & 5/12 &  3 + 2\sqrt{3} & 6 &  -3  \\
1/4 & 1/3 & \sqrt{3} & 0 & -3 \\
1/3 & 5/12 & (3+2\sqrt{3})/3 & 2 & -1/3 \\
\end{array}
$$
\centerline{Table 1: Quadratic ratios of tangents, that are non-units} 
\par

\section{Periodic points in genus two} \label{PPgenustwo}

The $\GL_2^+(\RR)$-action on $\Omega M_2$ respects the
multiplicities of the zeroes of the differential. We
denote the corresponding strata by $\Omega M_2(2)$
and $\Omega M_2(1,1)$.
\par
\subsection{The stratum $\Omega M_2(2)$}
\par
By \cite{McM04a} each Veech surface in the stratum
$\Omega M_2(2)$ contains an $L$-shaped surface in
its $\GL_2^+(\RR)$-orbit. 
\par
\centerline{\includegraphics{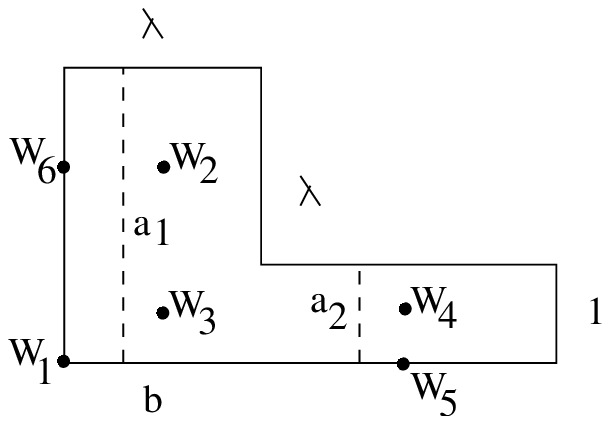}}
\centerline{Figure 1: 
Admissible representative of a primitive 
Veech surface in $\Omega M_2(2)$}
\par
Here $\lambda = (e + \sqrt(e^2 + 4b))/2$, where 
$e \in \{-1,0,1 \}$ and $b \in \NN$ with the 
restriction that $e+1 < b$ and if $e=1$ then $b$ is even. 
A triple $(b,e,\lambda)$ satisfying these conditions
is called {\em admissible}.
\par
In this section we show:
\par
\begin{Thm} \label{PerStratum2}
The only periodic points on a primitive Veech
surface in $\Omega M_2(2)$ are the Weierstra\ss\ points.
\end{Thm}
\par
{\bf Proof:} 
Let $f: X_1 \to C_1$ be a fibred surface (and $\ol{X_1} \to \ol{C_1}$
its stable completion) 
corresponding to a given primitive Veech surface in $\Omega M_2(2)$. 
Suppose the Veech surface contains a non-Weierstra\ss\ 
periodic point $P$. Passing to an unramified cover of $C_1$ we may suppose
by Lemma \ref{perfib} that the Weierstra\ss\ points and $P$
extend to sections $s_{W_i}$ for $i=1,\ldots,6$ and $s_P$ of $f$.
We suppose that $s_{W_1}$ passes through the zero of $\omega^0$.
By Thm.~\ref{peristor} and the primitivity
assumption the section $s_P - s_{W_1}$ is torsion.
\par
We start from an admissible representative $(X^0, \omega^0)$ 
and observe what happens if we degenerate
in the vertical direction, i.e.\ along the
path $$(X^t, \omega^t) = \diag(e^t, e^{-t})\cdot (X^0,\omega^0)$$ 
for $t \to \infty$ , where $(1\;\; 0)^T$
corresponds to the horizontal and $(0\;\;1)^T$ to the
vertical direction.
Note that the action of $\diag(e^t, e^{-t})$ 
does not change the ratio of the heights of the
vertical cylinders. 
\newline
The curves $X^t$ are obtained by first cutting $X^0$ along the
centers of the vertical cylinders and then glueing a pair of annuli (of
some modulus increasing with $t$) along the cut circles.
The limit curve is obtained by 'squeezing' the interior
of each vertical cylinder to a point. See \cite{Ma75} for more
details.
\par
Hence the stable model of the limit curve $X^\infty$
is a rational curve with two pairs of points identified.
By construction these nodes are fixed points of the
hyperelliptic involution.
\par
Normalising suitably, we may suppose that $X^\infty$ looks as follows: 
The Weierstra\ss\ section $s_{W_1}$ intersects $X^\infty$
in $\infty$, $s_{W_2}$ intersects $X^\infty$ in zero
and the hyperelliptic involution becomes
$z \mapsto -z$. By a linear transformation
we may suppose that $s_P$ intersects $X^\infty$ in $1$
and that the remaining Weierstra\ss\ sections
are glued to pairs $x$ with $-x$ and $y$ with $-y$.
for some $x,y \in \CC \sms \{0,\pm 1\}$.
\par
Furthermore $\omega^0$ comes from a subbundle of
$f_* \omega_{\ol{X_1}/\ol{C_1}}$ and specialises to
a section $\omega^\infty$ of the dualizing sheaf on the 
singular fibre $X^\infty$. Thus the differential
$\omega^\infty$ has to
vanish to the order two at $\infty$ and it  
has simple poles at $x$, $-x$, $y$ and $-y$. Furthermore
we must have 
$${\rm Res}_{\omega^\infty}(x) = {\rm Res}_{\omega^\infty}(-x) 
\quad \text{and} \quad
{\rm Res}_{\omega^\infty}(y) = {\rm Res}_{\omega^\infty}(-y). $$
The differential
$$ \omega^\infty = (\frac{y}{z-x} - \frac{y}{z+x} - \frac{x}{z-y}
+ \frac{x}{z+y})dz $$
has this property and by Riemann-Roch it is unique up 
to scalar multiple. 
\par
The invariance of the height ratios implies
$$ \frac{y}{x} = \frac{{\rm Res}_{\omega^\infty}(x)}
{{\rm Res}_{\omega^\infty}(-y)}
= \frac{\int_{a_1} \omega^0}{\int_{a_2} \omega^0} = \frac{\lambda+1}{1}$$
up to sign and interchanging the roles of $x$ and $y$. 
\par
Due to the irreducibility of $X^\infty$ and \cite{McM04c} Thm.\ 3.4 
the divisor $s_P - s_{W_1}$ remains torsion on the singular fibre.
Thus $1 -\infty$ is say $N$-torsion.
In order to have a well-defined map $g: X^\infty \to \PP^1$,
$z \mapsto (z-1)^N$ we must have both
$$ (x-1)^N = (-x-1)^N \quad \text{and} \quad (y-1)^N = (-y-1)^N.$$
This implies that $x,y \in i\RR$, in fact we must have $x = i\tan(A\pi/N)$ and
$y = i\tan(B\pi/N)$ for some $A,B \in \ZZ$. 
\newline
We now use the list of tangent ratios of the previous section:
If $y/x$ is a unit then 
$$|\Norm(\lambda+1)| = |e+1-b| = 1.$$
The only admissible triples that satisfy this condition are
$e=0$ and $b=2$ (which gives the octogon) and $e=-1$ and $b=1$
(which gives the pentagon). For these two cases the theorem
was proved in \cite{GuHuSc03} Examples 3 and 4. Alternatively
one might use that these Teichm\"uller curves pass through
the Riemann surfaces $y^2 = x (x^4 -1)$ resp.\ $y^2 = x^5 -1$, 
whose torsion points are known (\cite{BoGr00}) . Then
one can conclude as in the proof of Thm.\ \ref{PerPtDec}.
If $y/x$ is not a unit we can rule out each element of the
list of table 1 by the conditions that $\Norm(\lambda+1)$ has
to be an integer and ${\rm Trace}(\lambda+1) \in\{1,2,3\}$.
\hfill $\Box$
\par

\subsection{The stratum $\Omega M_2(1,1)$}
\par
This stratum contains only one primitive Teichm\"uller curve,
the one generated by the regular decagon (see 
\cite{McM04b} and \cite{McM04c}). The decagon
corresponds to the Riemann surface $y^2 = x^5 -1$
with the differential $xdx/y$ and its $\GL^+_2(\RR)$-orbit 
contains the following translation surface:
\par
\centerline{\includegraphics{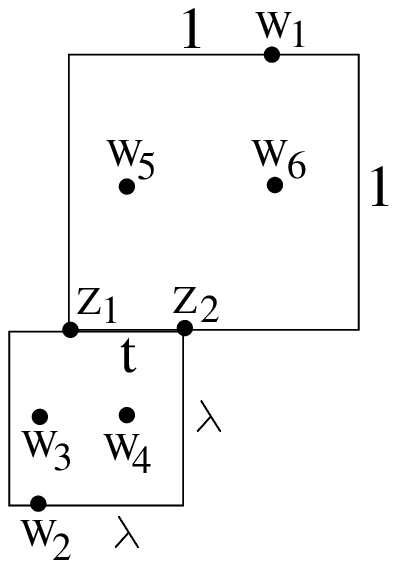}}
\centerline{Figure 2: Representative in the $\GL_2^+(\RR)$-orbit of 
the decagon}
Here the side lengths of the squares are $1$ and 
$\lambda = (-1+\sqrt{5})/2 $ and $t=\sqrt{5}/5$. We denote 
this surface by $(X^0, \omega^0)$.
\par
\begin{Thm} \label{PerPtDec}
The only periodic points of the decagon
are the Weierstra\ss\ points and the zeroes of
$\omega^0$.
\end{Thm}
\par
{\bf Proof:} We consider the degeneration of the
surface in Figure $2$ by squeezing the horizontal
direction. The singular fibre $X^\infty$ in the limit
is a $\PP^1$ with two pairs of points identified.
We suppose that the sections $s_{Z_i}$ extending the
two zeroes $Z_1$ and $Z_2$ of $\omega^0$ intersect
$X^\infty$ in $0$ and $\infty$ respectively. We 
may suppose that the hyperelliptic involution acts
by $z \mapsto 1/z$. This forces the Weierstra\ss\ 
sections $s_{W_1}$ and $s_{W_2}$ to intersect in $X^\infty$
in $-1$ and $1$ respectively. Furthermore $s_{W_3}$ and $s_{W_4}$ intersect
$X^\infty$ in the pair of identified points $x$ and $1/x$ 
while $s_{W_5}$ and $s_{W_6}$ degenerate to the
pair of identified points $y$ and $1/y$.
\par
The differential $\omega^\infty$ has simple zeros in $0$ and
$\infty$ and simple poles at $x$, $1/x$, $y$ and $1/y$ such that
$${\rm Res}_{\omega^\infty}(x) = {\rm Res}_{\omega^\infty}(1/x) 
\quad \text{and}\quad
{\rm Res}_{\omega^\infty}(y) = {\rm Res}_{\omega^\infty}(1/y). $$
This implies that 
$$ \omega^\infty = (\frac{y - 1/y}{z-x} - \frac{y - 1/y}{z - 1/x} + 
\frac{1/x - x}{z-y}
- \frac{1/x - x}{z - 1/y})dz $$
\newline
By Thm.\ \ref{peristor} the difference $Z_1 -Z_2$ is torsion. 
Considering the surface $y^2 =x^5 -1$ one notices (see e.g.\ \cite{BoGr00})
that $Z_1-Z_2$ is $5$-torsion, that $Z_1-W_1$ is $5$-torsion and
that whenever $R- W_1$ is torsion with $R$ a non-Weierstra\ss\ 
point, then $R - W_1$ is $5$-torsion.
Hence $g_1(z) = z^5$ and $g_2(z) = (z+1)^5$ are well-defined on $X^\infty$. 
This implies that (up to replacing $x$ or $y$ by its inverse or
interchanging $x$ and $y$) we
have $x = \exp(2\pi i /5)$ and $y = \exp(4\pi i /5)$. If there
was another periodic point $R$ on the decagon, which becomes $r$ 
on $X^\infty$ then $g_3(z) = (z-r)^5$ would have to be well-defined
on $X^\infty$. This implies immediately that $r$ has to be real. Let
$M_{10}$ be the set of complex numbers with argument a multiple
of $2\pi/10$. Then $\{x - r, r\in \RR \}$ intersects $M_{10}$
only for $r = -1$ or $r \geq 0$ and  $\{y - r, r\in \RR \}$
intersects $M_{10}$ only for $r \leq 0$. Hence there is no
possible choice for  $r \not\in \{-1,0,\infty\}$.
\hfill $\Box$
\par


\end{document}